\documentclass[11pt, reqno]{amsart}

 \usepackage{amsmath}
 \usepackage{amssymb}
\usepackage{bm}

\newcommand{\mysection}[1]{\section{#1}
      \setcounter{equation}{0}}

\newcommand{\nlimsup}{\operatornamewithlimits{\overline{lim}}}

\newtheorem{theorem}{Theorem}[section]
\newtheorem{lemma}[theorem]{Lemma}

\newtheorem{corollary}[theorem]{Corollary} 

\theoremstyle{definition}

\theoremstyle{remark}
\newtheorem{remark}{Remark}[section]

\newtheorem{example}{Example}[section]

\newcommand{\tr}{\text{\rm tr}\,}

 \makeatletter 
 \def\dashint{%
 \operatorname%
 {\,\,\text{\bf--}\kern-.98em\DOTSI\intop\ilimits@\!\!}}
\def\ninf{\qopname\relax\@empty{inf\phantom{p}\!\!\!}}
 \makeatother

\newcommand\bbeta{\text{\raise-.2ex\hbox{$\bm{\beta}$}}}

\newcommand\bR{\mathbb{R}}

\newcommand\bS{\mathbb{S}}

\newcommand\cF{\mathcal{F}}

\newcommand\cN{\mathcal{N}}
\newcommand\cS{\mathcal{S}}

\def\sft{{\sf t}}
\def\sfT{{\sf T}}

\begin{document}

\title[Stochastic It\^o equations with drift in $L_{d+1}$]
{On time inhomogeneous stochastic It\^o equations 
with drift in $L_{d+1}$}

\author{N.V. Krylov}
 
\email{nkrylov@umn.edu}
\address{127 Vincent Hall, University of Minnesota,
 Minneapolis, MN, 55455}

 \dedicatory{To the memory of A.V. Skorokhod}



\begin{abstract}
We prove the solvability of It\^o stochastic
equations with uniformly nondegenerate, bounded,
measurable diffusion and drift in $L_{d+1}(\bR^{d+1})$.
Actually, the powers of summability
of the drift in $x$ and $t$
could be different. Our results seem to be new even if
the diffusion  is constant. The method of proving
the solvability belongs to A.V. Skorokhod.
Weak uniqueness of solutions is an open problem even
if the diffusion is constant.

\end{abstract}

\maketitle

\mysection{Introduction}

Let $\bR^{d}$ be a  Euclidean space of points
$x=(x^{1},...,x^{d})$, $d\geq 2$. We fix some 
$p,q\in[1,\infty]$ such that  
\begin{equation}
                                                    \label{5.10.1}
  \frac{d}{p}+\frac{1}{q}\leq 1
\end{equation}
with further restrictions on them to be specified later.
The goal of this article is to study the solvability of
It\^o's stochastic equations of the form
\begin{equation}
                                                 \label{11.29.2}
x _{t}=x^{(0)}  +\int_{0}^{t}\sigma (t^{(0)}+s,x_{s})\,dw_{s}
+\int_{0}^{t}b (t^{(0)}+s,x_{s}) \,ds,
\end{equation}
where $w_{t}$ is a $d$-dimensional Wiener process,
$\sigma$ is a uniformly nondegenerate, bounded,
Borel function with values in the set of symmetric
$d\times d$ matrices, $b$ is a Borel measurable $\bR^{d}$-
valued function given on $(-\infty,\infty)\times\bR^{d}$
such that
\begin{equation}
                                                    \label{5.10.2}
\int_{\bR}\Big(\int_{\bR^{d}}|b(t,x)|^{p}\,dx\Big)^{q/p}\,dt<\infty 
\end{equation}
if $p\geq q$ or
$$
\int_ {\bR^{d}} \Big(\int_{\bR}|b(t,x)|^{q}\,dt\Big)^{p/q}\,dx<\infty 
$$
if $ p\leq q$. If $p=\infty$ or $q=\infty$ we interpret this conditions
in a natural way.
Observe that the case $p=q=d+1$ is not excluded and in this case
the condition becomes $b\in L_{d+1}(\bR^{d+1})$.
Under this condition the solvability of (1.2) was proved
in \cite{GM_01}.

We are talking, of course,
about weak solutions and prove their existence in Theorem
\ref{theorem 9.6.4}. In Theorem \ref{theorem 4.27.1}
we prove the existence of strong Markov processes
corresponding to diffusion $\sigma$ and drift $b$ 
with the above properties. 
If $b$ is bounded, as we know from \cite{Ya_20},  
there exist strong Markov and strong Feller processes
with diffusion $\sigma$ and drift $b$
for which the Harnack inequality holds
and the caloric functions are H\"older continuous.
We are far from proving such  fine  properties.

The main technical tools
are collected in Section \ref{section 5.10.1}
where we prove new mixed norms estimates
for the distributions of semimartingales.
The treatment there, actually, follows very closely
the work by A.I. Nazarov \cite{Na_15} written in terms
of PDEs.

There is a vast literature about stochastic equations with irregular
drift. Probably one of the first authors starting this area was
N.I. Portenko, see his book \cite{Po_82}, where he constructed
diffusion processes with sufficiently regular $\sigma$ and $b\in L_{p}
(\bR^{d+1})$, $p>d+2$. This condition on $b$ was later refined
in many articles with various ambitious 
goals in them to the requirement that
$b$ be such that \eqref{5.10.2} holds not
under condition \eqref{5.10.1}  but rather
\begin{equation}
                                                    \label{5.10.3}
\frac{d}{p}+\frac{2}{q}\leq 1.
\end{equation}
We refer the reader to the recent articles
\cite{Na_18},
\cite{BFGM_19},   \cite{XZ_20}  and the references therein
for the discussion of many powerful results obtained
 under condition \eqref{5.10.3}, when the case
of equality is treated as ``critical''.
It could be critical in some respects
but not for obtaining our results, that seem to be the first
ones about the existence of solutions and Markov
processes with our condition on the drift.
Still it is worth emphasizing that our
condition is different if $p\geq q$
(and hence $p\geq d+1$) or $p<q$, whereas there is no such
  distinction attached to 
  \eqref{5.10.3}.

 We assume that $d\geq2$ and denote
$$
B_{R} =\{x\in\bR^{d}:| x|<R\}, 	
\quad D_{i}=\frac{\partial}{\partial x^{i}},
\quad D_{ij}=D_{i}D_{j}\quad \partial_{t}=\frac{\partial}{\partial t}.
$$
For $p,q\in[1,\infty]$ we introduce the space $L_{p,q}$ as the space of Borel
functions on $\bR^{d+1}$ such that
$$
\|f\|^{q}_{p,q}:=
\int_{\bR}\Big(\int_{\bR^{d}}|f(t,x)|^{p}\,dx\Big)^{q/p}\,dt<\infty 
$$
if $p\geq q$ or
$$
\|f\|^{p}_{p,q}:=\int_ {\bR^{d}} \Big(\int_{\bR}|f(t,x)|^{q}\,dt\Big)^{p/q}\,dx<\infty 
$$
if $ p\leq q$ with natural interpretation
of these definitions if $p=\infty$ or $q=\infty$.
To better memorize these formulas observe that
$p $ is associated with integration with respect to $x$,
$q$ with that with respect to $t$ and the interior
integral is always elevated to the power $\leq 1$.
 In case $p=q=d+1$ we abbreviate $L_{d+1,d+1}=L_{d+1}$,
$\|\cdot\|_{d+1,d+1}=\|\cdot\|_{d+1}$.
 
\mysection{An example of nonexistence}
                                            \label{section 4.23.1}

\begin{example}
                                             \label{example 9.6.1}
Suppose that   numbers $\alpha$ and $\beta$ satisfy 
\begin{equation}
                                                  \label{4.18.1}
0<\alpha\leq \beta <1,\quad \alpha+\beta=1
\end{equation}
 and set
$$
b(t,x)=-\frac{1}{t^{\alpha}|x |^{\beta}}\frac{x }{|x |}
I_{0<|x|\leq 1,t\leq 1}.
$$
Observe that if $d/p+1/q=1+\varepsilon$, $\varepsilon>0$, 
one can take $\beta=d/(p+p\varepsilon)$, $\alpha=1/(q+q\varepsilon)$
and then 
$$
\int_{0}^{1}\Big(\int_{|x|\leq 1}|b(t,x)|^{p}\,dx\Big)^{q/p}dt<\infty,
\quad 
\int_{|x|\leq 1}\Big(\int_{0}^{1}|b(t,x)|^{q}\,dt\Big)^{p/q}dx<\infty.
$$
Also note
 that if  $p\leq qd$ (say $p=q$), condition \eqref{4.18.1}  is satisfied.

However, it turns out that no matter which $\alpha,\beta$
we take satisfying \eqref{4.18.1}    there is no solutions of the equation $dx_{t}=dw_{t}
+b(t,x_{t})\,dt$ starting at zero, where $w_{t}$ is a $d$-dimensional
Wiener process.

To prove this assume the contrary. Namely, assume the there is a stopping
time $\tau$ such that $P(\tau>0)>0$ and for $t\leq \tau$
there is $x_{t}$ such that
$$
x_{t}=w_{t}+\int_{0}^{t}b(s,x_{s})\,ds.
$$
We may assume that $\tau\leq 1$ and before $\tau$ the process 
is in $B_{1}$.
Then for $t\leq \tau$ 
\begin{equation}
                                                \label{4.18.4}
dx_{t}=-\frac{1}{t^{\alpha}|x_{t}|^{\beta}}\frac{x_{t}}{|x_{t}|}
I_{ x_{t}\ne 0 }\,dt+dw_{t},
\end{equation}
$$
d|x_{t}|^{2}=-2\frac{|x_{t}|}{t^{\alpha}|x_{t}|^{\beta}}\,dt
+d\,dt+2x_{t}\,dw_{t}.
$$
 
We will be  interested in $|x_{t}|^{1+\beta}=\xi_{t}^{(1+\beta)/2}$,
where $\xi_{t}=|x_{t}|^{2}$. 
By It\^o's formula for any $\varepsilon>0$ we have
$$
d(\xi_{t}+\varepsilon)^{(1+\beta)/2}=
\frac{1+\beta}{2}(\xi_{t}+\varepsilon)^{( \beta-1)/2}\,d\xi_{t}
+\frac{\beta^{2}-1}{8}(\xi_{t}+\varepsilon)^{( \beta-3)/2} 4|x_{t}|^{2}\,dt
$$
\begin{equation}
                                                           \label{3.18.1}
=I_{t}(\varepsilon)\,dt+J_{t}(\varepsilon)\,dt 
+(1+\beta)(\xi_{t}+\varepsilon)^{( \beta-1)/2}x_{t}\,dw_{t},
\end{equation}
where
$$
I_{t}(\varepsilon)=-(1+\beta)(\xi_{t}+\varepsilon)^{( \beta-1)/2}
\frac{|x_{t}|^{\alpha}}{t^{\alpha}} ,
$$

$$
J_{t}(\varepsilon) =
\frac{1+\beta}{2} \Big[ d+(\beta-1)
(\xi_{t}+\varepsilon)^{-1}  |x_{t}|^{2} \Big](\xi_{t}+\varepsilon)^{
(\beta-1)/2}
$$
Since $(\xi_{t}+\varepsilon)^{-\alpha/2}
 |x_{t}|^{\alpha}\uparrow I_{x_{t}\ne0}$ as $\varepsilon\downarrow0$,
by the dominated convergence theorem
$$
\int_{0}^{t}I_{s}(\varepsilon)\,ds\to
-(1+\beta)\int_{0}^{t}I_{x_{s}\ne0}
\frac{1}{s^{\alpha}}\,ds,
$$
which is finite. 

Furthermore,  since $|x_{s}|^{ \beta-1 }x_{s}$ is bounded
on each trajectory, by the dominated
convergence theorem
$$
\int_{0}^{t}|(\xi_{s}+\varepsilon)^{(\beta-1)/2}x_{s}
-|x_{s}|^{\beta-1 }x_{s}|^{2}\,ds\to0,
$$
and we conclude from \eqref{3.18.1} that for $t\leq\tau$
$$
|x_{t}|^{1+\beta}=-(1+\beta)\int_{0}^{t}I_{x_{s}\ne0}
\frac{1}{s^{\alpha}}\,ds
$$
\begin{equation}
                                                       \label{4.18.2}
+\lim_{\varepsilon\downarrow0}\int_{0}^{t}J_{s}(\varepsilon)\,ds
+(1+\beta)\int_{0}^{t}|x_{s}|^{\beta-1}x_{s}I_{x_{s}\ne0}\,dw_{s}
\end{equation}
and the above limit exists and is finite. Since
$2J_{s}(\varepsilon)\geq (\xi_{s}+\varepsilon)^{(\beta-1)/2}$,
it follows that
$$
\int_{0}^{t}|x_{s}|^{\beta-1 }\,ds=
\lim_{\varepsilon\downarrow0}\int_{0}^{t}
(\xi_{s}+\varepsilon)^{(\beta-1)/2}\,ds
$$
and the left-hand side is finite. 
In particular, 
\begin{equation}
                                                    \label{4.18.5}
\int_{0}^{\tau}I_{x_{s}=0}\,ds=0.
\end{equation}

Now by the dominated convergence theorem
\eqref{4.18.2} implies that  
$$
|x_{t}|^{1+\beta}=-(1+\beta)\int_{0}^{t} 
\frac{1}{s^{\alpha}}\,ds
$$
$$
+(1/2)(1+\beta)\int_{0}^{t}(d+\beta-1)|x_{s}|^{\beta-1 }\,ds
+(1+\beta)\int_{0}^{t}|x_{s}|^{\beta-1}x_{s} \,dw_{s}.
$$

Next, use $\alpha\leq \beta$ and H\"older's inequality to
conclude that
$$
\int_{0}^{t}|x_{s}|^{-\alpha}\,ds=
\int_{0}^{t}\Big(\frac{1}{s^{\alpha}|x_{s}|^{\beta}}
\Big)^{\alpha/\beta}s^{\alpha^{2}/\beta}\,ds
$$
 $$
\leq\Big(\int_{0}^{t}\frac{1}{s^{\alpha}
|x_{s}|^{\beta}}\,ds\Big)^{\alpha/\beta}
\Big(\int_{0}^{t}s^{\alpha^{2}/(\beta-\alpha)}
\,ds\Big)^{(\beta-\alpha)/\beta}.
$$
Since, $\alpha^{2}/(\beta-\alpha)+1=(\alpha^{2}+1-2\alpha)/(\beta-\alpha)
=\beta^{2}/(\beta-\alpha)$
$$
\int_{0}^{t}|x_{s}|^{-\alpha}\,ds\leq N 
\Big(\int_{0}^{t}\frac{1}{s^{\alpha}
|x_{s}|^{\beta}}\,ds\Big)^{\alpha/\beta}t^{\beta},
$$
where $N =N (\alpha,\beta)$
(which is trivial if $\alpha=\beta$).
Thus,
$$
|x_{t}|^{1+\beta}+ct^{\beta}\leq N_{1}
\Big(\int_{0}^{t}\frac{1}{s^{\alpha}
|x_{s}|^{\beta}}\,ds\Big)^{\alpha/\beta}t^{\beta}+
(1+\beta)\int_{0}^{t}|x_{s}|^{\beta-1}x_{s}\,dw_{s},
$$
where $c>0$ is a constant.
For   equation \eqref{4.18.4} to make sense we should have
\begin{equation}
                                            \label{4.18.3}
\int_{0}^{\tau}\frac{1}{s^{\alpha}
|x_{s}|^{\beta}}\,ds<\infty 
\end{equation}
  (a.s.). Therefore 
$$
\gamma:=\tau\wedge\inf\{t\geq0:N_{1}
\Big(\int_{0}^{t}\frac{1}{s^{\alpha}
|x_{s}|^{\beta}}\,ds\Big)^{\alpha/\beta}\geq c/2\},
$$
is a stopping time such that $P(\gamma>0)=P(\tau>0)$.
It follows that  for any $t>0$
$$
\int_{0}^{t}I_{s<\gamma}|x_{s}|^{\beta-1}x_{s}\,dw_{s} 
 \geq 0,
$$
which is only possible if $I_{s<\gamma}|x_{s}|^{\beta-1}x_{s}=0$
for almost all $(\omega,s)$. Then $x_{s}=0$ for $s<\gamma$
and \eqref{4.18.5} is only possible if $P(\tau=0)=1$.
\end{example}

\mysection{An existence theorem}
                                            \label{section 4.23.2}

In this section  we state a result saying that
in a wide class of cases there exists a 
probability space and a Wiener process
on this space such that a stochastic 
equation having measurable coefficients as
well as this Wiener process is solvable. In other words,  
 according to conventional terminology, we are talking
 here about ``weak'' solutions
of a stochastic equation. The main difference between ``weak'' solutions and
usual (``strong'') solutions consists in
 the fact that the latter can be constructed
on any a priori given probability space on the basis of any given
Wiener process.

Let $\sigma(t,x)$ be Borel $d\times d$ symmetric
matrix valued, $b(t,x)$ be Borel $\bR^{d}$-valued functions
given on $\bR^{d+1} :=(-\infty,\infty)\times \bR^{d}$. 
We assume that the eigenvalues of $\sigma(t,x)$ are between
$\delta$ and $\delta^{-1}$, where $\delta\in(0,1)$ is a fixed number.
The set of such matrices we denote by $\bS_{\delta}$.

Next, fix   numbers $p,q\in(1,\infty)$, $\|b\|\in(0,\infty)$ and
let   $b^{n}(t,x)$, $n= 1,2,...$, be   $\bR^{d}$-valued 
Borel functions
on $\bR^{d+1}_{+}$ and suppose that     
$$
\|b \|_{p,q},
\|b^{n}\|_{p,q}\leq \|b\|,
\quad n= 1,2,...,\quad \frac{d}{p}+\frac{1}{q}=1
$$
and $b^{n}\to b $ as $n\to\infty$ in $L_{p,q} $.
 Let   $\sigma^{n}(t,x)$, $n= 1,2,...$,    
be  Borel functions on $\bR^{d}$ with values in $\bS_{\delta}$
  such that
$\sigma^{n}\to \sigma $ as $n\to\infty$ ($\bR^{d+1} $-a.e.).

\begin{theorem}
                                           \label{theorem 9.6.4}
Take $(t^{0},x^{0})\in\bR^{d+1}$. 
(i) There exists 
a probability space $(\Omega ,\cF ,P )$,
a filtration of $\sigma$-fields $\cF _{t}\subset \cF $, $t\geq0$,
a process $w _{t}$, $t\geq0$, which is a $d$-dimensional Wiener process
relative to $\{\cF _{t}\}$, and an $\cF _{t}$-adapted
process $x_{t}$ such that 
 (a.s.) for all   $t\geq0$ equation \eqref{11.29.2} holds.

(ii) Furthermore,
let $(t^{n},x^{n})\in\bR^{d+1} $, $n= 1,2,...$, and let $(t^{n},x^{n})
\to (t^{0},x^{0})$  as $n\to\infty$. 
Assume that for each $n=1,2,...$  
 there exists 
a probability space $(\Omega^{n},\cF^{n}, P^{n})$,
a filtration of $\sigma$-fields $\cF^{n}_{t}\subset \cF^{n}$, $t\geq0$,
a process $w^{n}_{t}$, $t\geq0$, which is a $d$-dimensional Wiener process
relative to $\{\cF^{n}_{t}\}$, and an $\cF^{n}_{t}$-adapted
process $x^{n}_{t}$ such that (a.s.) for all   $t\geq0$
\begin{equation}
                                                 \label{11.29.1}
x^{n}_{t}=x^{n} +\int_{0}^{t}\sigma^{n}(t^{n}+s,x^{n}_{s}) \,dw^{n}_{s}
+\int_{0}^{t}b^{n}(t^{n}+s,x^{n}_{s}) \,ds.
\end{equation}

Then  the finite dimensional distributions of  
a subsequence of 
$ x^{n}_{\cdot} $ converge  weakly to 
the corresponding distributions of one of the solutions of
\eqref{11.29.2}  described in (i). Moreover, if $p\geq q$,  
the set of distributions of $ x^{n}_{\cdot} $ on 
$C([0,\infty),\bR^{d })$ is tight.
\end{theorem}

The proof of this theorem, following a similar proof 
by A.V. Skorokhod, is given in Section
\ref{section 5.13.1}, after we make
a crucial step  
in the next section where we prove, in particular, that
for solutions of \eqref{11.29.2}, any
Borel $f\geq0$, and $T\in(0,\infty)$
\begin{equation}
                                                    \label{4.27.3}
E\int_{0}^{T}f(t,x_{t})\,dt\leq N\|f\|_{p,q},
\end{equation}
where $N$ is independent of $f$ and $(t^{0},x^{0})$.

It is worth saying that deciding whether the
solutions
of \eqref{11.29.2}
are weakly unique or not   under our conditions is a very
challenging open problem even if $\sigma^{ij}\equiv \delta^{ij}$.

 \begin{remark}
Theorem \ref{theorem 9.6.4} is also true if $d/p+1/q<1$.
This can  be seen from its proof which becomes somewhat
more technical in that case because of the form of our main estimate
\eqref{5.6.40}. Also the main interest in Theorem \ref{theorem 9.6.4}
is, of course, the lowest {\em local\/} integrability of $b$,
when the condition $d/p+1/q=1$ is weaker than $d/p+1/q<1$
due to H\"older's inequality.
\end{remark}

\mysection{Estimates of the distributions of semimartingales}
                                                \label{section 5.10.1}

Here we first prove a version of Lemma 5.1 of \cite{Kr_86}.
The proof given in \cite{Kr_86} uses somewhat advanced
knowledge of very powerful results from
the theory of fully nonlinear parabolic
equations. We give a proof based on a simpler fact
which in turn was one of the cornerstones of that theory.

Let $(\Omega,\cF,P)$ be a complete probability
space, let $\cF_{t}, t\geq0$, be an increasing family of
complete $\sigma$-fields $\cF_{t}\subset\cF$, $t\geq0$,
let $m_{t}$ be an $\bR^{d}$-valued continuous local martingale
relative to $\cF_{t}$, let  $A_{t}$ be a continuous
$\cF_{t}$-adapted
nondecreasing process, let $B_{t}$ be a continuous $\bR^{d}$-valued
$\cF_{t}$-adapted process which has finite variation (a.e.)
on each finite time interval. Assume that
$$
A_{0}=0,\quad m_{0}=B_{0}=0,\quad d\langle m
 \rangle_{t} \ll  dA_{t}
$$
and that we are also given progressively measurable relative to
$\cF_{t}$ nonnegative processes $r_{t}$ and $c_{t}$. Finally, take
an $\cF_{0}$ measurable $\bR^{d}$-valued $x_{0}$ and
introduce
$$
x_{t}=x_{0}+m_{t}+B_{t},\quad \tau_{t}=\int_{0}^{t}r_{s}\,dA_{s},
\quad \phi_{t}=\int_{0}^{t}c_{s}\,dA_{s},\quad
a^{ij}_{t}=\frac{1}{2}\frac{d\langle m^{i},m^{j}\rangle
_{ t}} 
{dA_{t}}.
$$
 
\begin{lemma}
                                            \label{lemma 4.23.3}
Let $\gamma$ be an $\cF_{t}$-stopping time and set
$$
A=E\int_{0}^{\gamma}e^{-\phi_{t}}\tr a_{s}\,dA_{t},\quad 
B=E\int_{0}^{\gamma}e^{-\phi_{t}}\,|dB_{t}|.
$$
Then for any Borel $f(t,x)\geq0$ we have
$$
E\int_{0}^{\gamma}e^{-\phi_{t}}(r_{t}\det a_{t})^{1/(d+1)}
f(\tau_{t},x_{t})\,dA_{t}
$$
\begin{equation}
                                            \label{4.24.1}
\leq N(d)(B^{2}+A)^{d/(2d+2)}\|f\|_{ d+1}.  
\end{equation}
 
\end{lemma}

Proof. Without losing generality we may assume that $A<\infty$
and $B<\infty$. Furthermore, just stopping the processes $A_{t}$,
$m_{t}$, and $B_{t}$ at time $\gamma$, we reduce the general case
to the one in which $\gamma=\infty$. In that case we also observe that,
as usual, it suffices to prove \eqref{4.24.1} for $f\in
C^{\infty}_{0}(\bR^{d+1})$.

After these reductions we use Theorem 2.2.4 of \cite{Kr_18}
according to which, for any $\lambda>0$
on $\bR^{d+1}$,
there exists a nonnegative function $v(t,x)$ such that  

(i)  all Sobolev derivatives  
 $\partial_{t}v $, $D_{i}v$, $ D_{ij}v$
exist and are bounded on $\bR^{d+1}$ and $v\le Ne^{- |x|/N}$
for all $t$, $x$ and a constant $N$; 

(ii) for any nonnegative symmetric $d\times d$ matrix
$\alpha $ and $r\ge0$,  
$$
r\partial_{t}v+\alpha^{ij}D_{ij}v -\lambda
(r+\tr\alpha)v+\sqrt[d+1]{r\det\alpha}\,f\le0,
 $$
\begin{equation}
                                                \label{eq:3.2.14}
\partial_{t}v-\lambda v\le0,\quad(\lambda 
v\delta_{ij}-D_{ij}v )\ge0,
\quad|Dv |\le\sqrt{\lambda}v\quad(a.e.) ,
\end{equation} 

(iii) for any $y\in\bR^{d}$, $t\in(-\infty,\infty)$,  
we have
\begin{equation}
                                                   \label{4.24.3}
v(t,  y)e^{-\lambda t}\le 
\,N(d)\frac{1}{\lambda^{d/(2d+2)}}I_{t},
\end{equation}
where
$$
I^{d+1}_{t}:= \int_{0}^{\infty}ds
\int_{\bR^{d}}e^{-\lambda(d+1)(t+s)}f^{d+1}(t+s, 
x)\,dx .
$$

Take a nonnegative $\zeta\in C^{\infty}_{0}(\bR^{d+1})$ with unit integral,
for $\varepsilon>0$ denote $\zeta_{\varepsilon}(t,x)=\varepsilon^{-(d+1)}
\zeta(\varepsilon t,\varepsilon x)$ and use the notation
$u^{(\varepsilon)}=u*\zeta_{\varepsilon}$. Then $v^{(\varepsilon)}$
is infinitely differentiable and in light of \eqref{eq:3.2.14},
for any nonnegative symmetric $d\times d$ matrix
$\alpha $ and $r\ge0$,  
$$
r\partial_{t}v^{(\varepsilon)}+\alpha^{ij}D_{ij}v^{(\varepsilon)} -\lambda
(r+\tr\alpha)v^{(\varepsilon)}
+\sqrt[d+1]{r\det\alpha}\,f^{(\varepsilon)}\le0,
 $$
\begin{equation}
                                                \label{4.25.1}
\partial_{t}v^{(\varepsilon)}-\lambda v^{(\varepsilon)}
\le0,\quad(\lambda 
v^{(\varepsilon)}\delta_{ij}-D_{ij}v^{(\varepsilon)} )\ge0,
\quad|Dv^{(\varepsilon)} |\le\sqrt{\lambda}v^{(\varepsilon)}  .
\end{equation} 
 Next, by It\^o's formula   the process
$$
v^{(\varepsilon)}(\tau_{t},x_{t})e^{-\phi_{t}-\lambda \tau_{t}}
-\int_{0}^{t}e^{-\phi_{s}-\lambda \tau_{s}}D_{i}v^{(\varepsilon)}
(\tau_{s},x_{s})\,dB^{i}_{s}
$$
$$
+\int_{0}^{t}e^{-\phi_{s}-\lambda \tau_{s}}\big((\lambda
r_{s}+c_{s})v^{(\varepsilon)}
-r_{s}\partial_{t}v^{(\varepsilon)}-a^{ij}_{s}
D_{ij}v^{(\varepsilon)}\big)(\tau_{s},x_{s})\,dA_{s}
$$
is a local martingale. Here owing to \eqref{4.25.1}
$$
((\lambda r_{s}+c_{s})v^{(\varepsilon)}
-r_{s}\partial_{t}v^{(\varepsilon)}-a^{ij}_{s}
D_{ij}v^{(\varepsilon)}\big)\,dA_{s}
-D_{i}v^{(\varepsilon)}\,dB^{i}_{s}
$$
$$
\geq (r_{s}\det a_{s})^{1/(d+1)}f^{(\varepsilon)}\,dA_{s}
-\lambda \tr a_{s}v^{(\varepsilon)}\,dA_{s}-
\sqrt{\lambda}v^{(\varepsilon)}\,|dB_{s}|.
$$
Therefore, for 
$$
M^{\varepsilon}=\sup_{t\geq 0,x\in\bR^{d}}
v^{(\varepsilon)}(t,x)e^{-\lambda t}
$$
 the process
$$
\kappa^{\varepsilon}_{t}:=
v^{(\varepsilon)}(\tau_{t},x_{t})
 e^{-\phi_{t}-\lambda \tau_{t}}+
\int_{0}^{t}e^{-\phi_{s}-\lambda \tau_{s}}
 (r_{s}\det a_{s})^{1/(d+1)} f^{(\varepsilon)}
(\tau_{s},x_{s})\,dA_{s}
$$
$$
-\int_{0}^{t}e^{-\phi_{s} }\big(
\lambda \tr a_{s} \,dA_{s} + 
\sqrt{\lambda} \,|dB_{s}|\big) \,M^{\varepsilon}
$$
 is a local supermartingale. In addition,
it is bounded from below by a summable quantity
($A,B<\infty$). Hence, it is a supermartingale
and by Fatou's lemma
$$
Ev^{(\varepsilon)}(0,x_{0})=\kappa^{\varepsilon}_{0}
\geq E\int_{0}^{\infty}e^{-\phi_{t}-\lambda \tau_{s}}
(r_{t}\det a_{t})^{1/(d+1)}
f^{(\varepsilon)}(\tau_{t},x_{t})\,dA_{t}
$$
$$
-M^{\varepsilon}(\lambda
A+\sqrt{\lambda}B).
$$
By sending $\varepsilon\downarrow 0$ and using \eqref{4.24.3} and
Fatou's lemma once more we obtain
$$
E\int_{0}^{\infty}e^{-\phi_{t}-\lambda \tau_{s}}
(r_{t}\det a_{t})^{1/(d+1)}
f (\tau_{t},x_{t})\,dA_{t}
$$
$$
\leq N(d)\frac{1}{\lambda^{d/(2d+2)}}\Big( 1
+\lambda
A+\sqrt{\lambda}B\Big)I_{0}.
$$
We replace here $e^{-\lambda t}f$ by $f$ and arrive at
$$
E\int_{0}^{\infty}e^{-\phi_{t} }
(r_{t}\det a_{t})^{1/(d+1)}
f (\tau_{t},x_{t})\,dA_{t}
$$
$$ \leq
N(d)\frac{1}{\lambda^{d/(2d+2)}} \Big(1
+\lambda
A+\sqrt{\lambda}B\Big)\|f\|_{d+1}.
$$
Now we use the arbitrariness of $\lambda$. If $A< B^{2}$,
then for $\lambda=B^{-2}$ we have
$$
\frac{1}{\lambda^{d/(2d+2)}}\Big(1 
+\lambda
A+\sqrt{\lambda}B\Big)\leq 3B^{d/(d+1)}\leq 3(B^{2}+A)^{d/(2d+2)}.
$$
If $A\geq B^{2}$ and $A>0$, then for $\lambda=A^{-1}$
the above inequality between the extreme terms still holds.
Finally, if $A=B=0$, then the left-hand side of
\eqref{4.24.1} is zero. The lemma is proved.

\begin{lemma}
                                            \label{lemma 5.6.1}
In the notation of Lemma \ref{lemma 4.23.3}
for any Borel $f( x)\geq0$ we have
\begin{equation}
                                            \label{5.6.1}
E\int_{0}^{\gamma}e^{-\phi_{t}} (\det a_{t}) ^{1/d}
f( x_{t})\,dA_{t}\leq N(d)(B^{2}+A)^{1/2}\|f\|_{L_{d }(\bR^{d })}.
\end{equation}
\end{lemma}

Proof. We follow a probabilistic version of
an argument in \cite{Na_15}. We again may concentrate on
the case of $A+B<\infty$, $\gamma=\infty$,
 and $f\in C^{\infty}_{0}(\bR^{d})$.
In that case observe that by Theorem 2.2.3 of \cite{Kr_18}
there exists a nonnegative function
$v(x)$ defined on $\bR^{d}$ such that 

(a) $v\le Ne^{- |x|/N}$ for all $x$ and a constant
$N$; the generalized derivatives   $D_{i}v
$ and $D_{ij}v $,
$i,j=1,..., d$,  are bounded on $\bR^{d}$;

(b) for any nonnegative symmetric $d\times d$ matrix
$\alpha$ (a.e.)  
$$
-\lambda v\,\tr\,\alpha+\alpha^{ij}D_{ij}v +\sqrt[d]{\det\alpha}f\le0,
\quad|Dv|\le\sqrt{\lambda}v,
 $$
\begin{equation}
                                                     \label{eq:3.2.9}
(\lambda v\delta_{ij}-D_{ij}v )\ge0;
\end{equation}

(c) for any $x\in\bR^{d}$  
\begin{equation}
                                                      \label{eq:3.2.10}
 v(x)   \le N(d)\lambda^{- 1/2}\Vert 
f\Vert _{L_{d}(\bR^{d})}.
\end{equation}

Then we closely follow the proof of Lemma \ref{lemma 4.23.3}.
Take a nonnegative $\zeta\in C^{\infty}_{0}(\bR^{d })$ with unit integral,
for $\varepsilon>0$ denote $\zeta_{\varepsilon}( x)=\varepsilon^{-d }
\zeta( \varepsilon x)$ and use the notation
$u^{(\varepsilon)}=u*\zeta_{\varepsilon}$. Then $v^{(\varepsilon)}$
is infinitely differentiable and in light of \eqref{eq:3.2.9},
for any nonnegative symmetric $d\times d$ matrix
$\alpha $,  
$$
 \alpha^{ij}D_{ij}v^{(\varepsilon)} -\lambda
 \tr\alpha v^{(\varepsilon)}
+\sqrt[d ]{ \det\alpha}\,f^{(\varepsilon)}\le0,
 $$
\begin{equation}
                                                \label{4.25.10}
 (\lambda 
v^{(\varepsilon)}\delta_{ij}-D_{ij}v^{(\varepsilon)} )\ge0, 
\quad|Dv^{(\varepsilon)} |\le\sqrt{\lambda}v^{(\varepsilon)}  .
\end{equation} 
 Next, by It\^o's formula   the process
$$
v^{(\varepsilon)}( x_{t})e^{-\phi_{t} }
-\int_{0}^{t}e^{-\phi_{s} }D_{i}v^{(\varepsilon)}
( x_{s})\,dB^{i}_{s}
$$
$$
+\int_{0}^{t}e^{-\phi_{s}-\lambda \tau_{s}}\big( 
 c_{s} v^{(\varepsilon)}
 -a^{ij}_{s}
D_{ij}v^{(\varepsilon)}\big)( x_{s})\,dA_{s}
$$
is a local martingale. Here owing to \eqref{4.25.10}
$$
\big( 
 c_{s} v^{(\varepsilon)}
 -a^{ij}_{s}
D_{ij}v^{(\varepsilon)}\big) \,dA_{s}
-D_{i}v^{(\varepsilon)}\,dB^{i}_{s}
$$
$$
\geq ( \det a_{s})^{1/ d }f^{(\varepsilon)}\,dA_{s}
-\lambda \tr a_{s}v^{(\varepsilon)}\,dA_{s}-
\sqrt{\lambda}v^{(\varepsilon)}\,|dB_{s}|.
$$
Therefore, for 
$$
M^{\varepsilon}=\sup_{ x\in\bR^{d}}
v^{(\varepsilon)}( x) 
$$
 the process
$$
\kappa^{\varepsilon}_{t}:=
v^{(\varepsilon)}( x_{t})
 e^{-\phi_{t} }+
\int_{0}^{t}e^{-\phi_{s} }f^{(\varepsilon)}
( x_{s})\,dA_{s}
$$
$$
-\int_{0}^{t}e^{-\phi_{s} }\big(
\lambda \tr a_{s} \,dA_{s}-
\sqrt{\lambda} \,|dB_{s}|\big) \,M^{\varepsilon}
$$
 is a local supermartingale. In addition,
it is bounded from below by a summable quantity
($A,B<\infty$). Hence, it is a supermartingale
and by Fatou's lemma
$$
Ev^{(\varepsilon)}( x_{0})=\kappa^{\varepsilon}_{0}
\geq E\int_{0}^{\infty}e^{-\phi_{t} }
( \det a_{t})^{1/d}
f^{(\varepsilon)}( x_{t})\,dA_{t}
$$
$$
-M^{\varepsilon}(\lambda
A+\sqrt{\lambda}B).
$$
By sending $\varepsilon\downarrow 0$ and using \eqref{eq:3.2.10} and
Fatou's lemma once more we obtain
$$
E\int_{0}^{\infty}e^{-\phi_{t} }
( \det a_{t})^{1/d}
f ( x_{t})\,dA_{t}
$$
$$
\leq N(d)\frac{1}{\lambda^{1/2}} \Big(1
+\lambda
A+\sqrt{\lambda}B\Big)\Vert 
f\Vert _{L_{d}(\bR^{d})}.
$$
 
Now we use the arbitrariness of $\lambda$. If $A< B^{2}$,
then for $\lambda=B^{-2}$ we have
$$
\frac{1}{\lambda^{1/2}} \Big(1
+\lambda
A+\sqrt{\lambda}B\Big)\leq 3B^{1/2}\leq 3(B^{2}+A)^{1/2}.
$$
If $A\geq B^{2}$ and $A>0$, then for $\lambda=A^{-1}$
the above inequality between the extreme terms still holds.
Finally, if $A=B=0$, then the left-hand side of
\eqref{5.6.1} is zero. The lemma is proved.

\begin{theorem}
                                         \label{theorem 5.5.1}
Assume  the notation of Lemma \ref{lemma 4.23.3} and let
 $p,q\in[1,\infty]$ be
such that
$$
 \theta:=1-\frac{d}{p}-\frac{1}{q}\geq0
$$
then for any Borel $f(t, x)\geq0$   we have
\begin{equation}
                                                   \label{5.6.40}
I(p,q,f):= E\int_{0}^{\gamma}e^{-\phi_{t}}\kappa_{t}
f(r_{t},x_{t})\,dA_{t}\leq N(d) (A+B^{2})^{ d/(2p)}
\|f\|_{p,q},
\end{equation}
where 
$\kappa_{t}=
r_{t}^{1/q}(\det a_{t})^{1/p}c^{\theta}_{t}$ 
and for any $\alpha\geq0$ we set $\alpha^{0}=1$ (say, if $\theta=0$).
\end{theorem}

Proof. By H\"older's inequality, if $\theta>0$,
$$
I(p,q,f)\leq  \Big(I(p (1-\theta ),q (1-\theta ),
f^{1/(1-\theta)})
\Big)^{1-\theta}.
$$
It follows that it suffices to concentrate on   $\theta=0$.
Then we observe that if $q=\infty$, then $p=d$ and
$$
\|f\|^{p}_{p,q}=\int_{\bR^{d}}\sup_{t\geq 0}f^{d}(t,x)\,dx.
$$
In that case \eqref{5.6.4} follows from Lemma \ref{lemma 5.6.1}.  
If $p=\infty$, then $q=1$,
and
$$
I(p,q,f)= E\int_{0}^{\gamma}r_{t} f (\tau_{t},x_{t})
\,dA_{t}\leq E\int_{0}^{\gamma} \sup_{x}f (\tau_{t},x)
\,d\tau_{t}
$$
$$
\leq \int_{0}^{\infty}\sup_{x}f (t,x)\,dt=\|f\| _{p,q}.
$$
In the third simple situation when $q=p=d+1$ estimate
\eqref{5.6.4} follows from Lemma \ref{lemma 4.23.3}.
We prove the lemma in the remaining cases by interpolating between
the above ones.

If  $p>q$ (and hence $p>d+1$) we take a
  nonnegative
function $h(t )$ such that $(hf)/h=f$ ($0/0:=0$) and
 use
$$
 r_{t}^{1/q}(\det a_{t})^{1/p}f=
 \Big(r_{t}^{1/q-1/p}h^{-1}\Big)\Big((r\det a_{t})^{1/p}fh\Big)
$$
along with H\"older's inequality. By  
performing simple manipulations we find
$$
I(p,q,f)\leq IJ
$$
\begin{equation}
                                                   \label{5.6.50}
:= 
\Big(I(\infty,1,
h^{-p /(p-d-1)})
\Big)^{ (p-d-1)/p}\Big(I(d+1 ,d+1  ,
(hf)^{p/(d+1)})
\Big)^{(d+1)/p}.
\end{equation}
Here  
$$
I\leq \Big(\int_{0}^{\infty}h^{-p /(p-d-1)}(t)\,dt\Big)^{ (p-d-1)/p}.
$$
Also
$$
J\leq
N(d)(B^{2}+A)^{d/(2p)}\|(hf)^{p/(d+1)}\|^{(d+1)/p}_{d+1}
$$
$$
=N(d)(B^{2}+A)^{d/(2p)}\Big(\int_{0}^{\infty} 
\Big(\int_{\bR^{d}}f^{p}(t,x)\,dx\Big)h^{p}(t)\,dt\Big)^{1/p}.
$$
We now choose $h$ so that
$$
h^{-p /(p-d-1)}(t)=\Big(\int_{\bR^{d}}f^{p}(t,x)\,dx\Big)h^{p}(t).
$$
Then both quantities become
$$
\Big(\int_{\bR^{d}}f^{p}(t,x)\,dx\Big)^{q/p},\quad
J\leq N(d)(B^{2}+A)^{d/(2p)}\|f\|^{q/p}_{p,q},\quad
I\leq \|f\|^{q(p-d-1)/p}_{p,q}
$$
and coming back to \eqref{5.6.50} we get \eqref{5.6.40}.

In the remaining case $q>p$ (and $q>d+1$)
we use
$$
r_{t}^{1/q}(\det a_{t})^{1/p}f=
 \Big((\det a_{t})^{1/p-1/q}h^{-1}\Big)\Big((r\det a_{t})^{1/q}fh\Big).
$$
This time for $h=h(x)$
$$
I(p,q,f)\leq IJ
$$
\begin{equation}
                                                   \label{5.6.6}
:= 
\Big(I(d ,\infty,
h^{-q /(q-d-1)})
\Big)^{ (q-d-1)/q}\Big(I(d+1 ,d+1  ,
(hf)^{q/(d+1)})
\Big)^{(d+1)/q}.
\end{equation}
 Here
$$
I\leq N(d)(B^{2}+A)^{(d/p-d/q)(1/2)}\Big(\int_{\bR^{d}}h^{-qd/(q-d-1)}
(x)\,dx\Big)^{(q-d-1)/(qd )},
$$
$$
J\leq N(d)(B^{2}+A)^{d/(2q) }\Big(\int_{\bR^{d}}h^{q}(x)
\Big(\int_{0}^{\infty}f^{q}(t,x)\,dt\Big)\,dx\Big)^{1/q}.
$$
We choose $h$ so that
$$
h^{-qd/(q-d-1)}
(x)=h^{q}(x)
\Big(\int_{0}^{\infty}f^{q}(t,x)\,dt\Big)
$$
and then easily come to \eqref{5.6.4}. The theorem
is proved.

\begin{corollary}
                                         \label{corollary 5.1.1}
Introduce a measure (Green's measure) on Borel 
subsets $\Gamma$ of $\bR^{d+1}$
by the formula
$$
G(\Gamma)=
E\int_{0}^{\gamma}e^{-\phi_{t}}\kappa_{t}
I_{\Gamma}(\tau_{t},x_{t})\,dA_{t}.
$$
Assume that $A,B<\infty$ and set $p'=p/(p-1),q'=q/(q-1)$.
 Then $G(\Gamma)$ is absolutely continuous and
its  density $G(t,x)$ is such that, if $p\geq q$,
$$
\Big(\int_{0}^{\infty}\Big(\int_{\bR^{d}}G^{p'}(t,x)\,dx\Big)^{q'/p'}
dt\Big)^{1/q'}
$$
and, if $p\leq q$,
$$
\Big(\int_{\bR^{d}}\Big(\int_{0}^{\infty}G^{q'}(t,x)\,dt\Big)^{p'/q'}
dx\Big)^{1/q'}
$$
is dominated by
$$
  N(d)(B^{2}+A)^{(1-\theta)d/(2p)}.
$$
\end{corollary}

\begin{theorem}
                                         \label{theorem 5.5.2}
 Under the assumptions of Theorem \ref{theorem 5.5.1}
let $p_{0} \in[1,\infty]$ and $q_{0}\in[1,\infty)$ be
such that
$$
\theta_{0}:=1-\frac{d}{p_{0}}-\frac{1}{q_{0}}\geq0.
$$
Also assume  that $d|B_{t}|\ll dA_{t}$
 and there exists a Borel $h(t,x)$
such that   ($P(d\omega)\times dA_{t}$-a.e.)
$$
|b_{t}|\leq \kappa^{0}_{t}h(\tau_{t}
,x_{t}),
$$
where $b_{t}=dB_{t}/dA_{t}$ and $\kappa^{0}_{t}=
r_{t}^{1/q_{0}}(\det a_{t})^{1/p_{0}}c^{\theta_{0}}_{t}$. 
Then for any Borel $f(t, x)\geq0$   we have
\begin{equation}
                                                   \label{5.6.4}
I(p,q,f):= E\int_{0}^{\gamma}e^{-\phi_{t}}\kappa_{t}
f( \tau _{t},x_{t})\,dA_{t}\leq N(d,p_{0},q_{0})C
\|f\|_{p,q},
\end{equation}
where 
$$
\kappa_{t}=r_{t}^{1/q}(\det a_{t})^{1/p}c^{\theta}_{t},\quad
C= \Big(A+ \|h\|_{p_{0},q_{0}}
^{ 2p_{0}/(p_{0}-d)   }
\Big)^{ d/(2p)} 
$$ and
  for any number $\alpha\geq0$ we set $\alpha^{0}=1$ (say, if $\theta=0$).

\end{theorem}

Proof. Observe that $p_{0}>d$ since $q_{0}<\infty$.
Then, we may assume that $A<\infty$ and $\|h\|_{p_{0},q_{0}}<\infty$.
Using stopping times we easily reduce the general situation
to the one in which $B<\infty$. After that, in light of 
Theorem
\ref{theorem 5.5.1}, we  need only prove that
\begin{equation}
                                                   \label{5.6.5}
B\leq N(d,p_{0},q_{0})\Big(A^{1/2}+ \|h\|_{p_{0},q_{0}}
 ^{ p_{0}/(p_{0}-d)  }\Big).
\end{equation}

By Theorem
\ref{theorem 5.5.1}
$$
B=E\int_{0}^{\tau}e^{-\phi_{t}}|dB_{t}|\leq
I(p_{0},q_{0},h)\leq N(d)
 (A+B^{2})^{ d/(2p_{0})}\|h\|_{p_{0},q_{0}}.
$$

Here if $B^{2}\leq A$, estimate \eqref{5.6.5} holds.
If $A\leq B^{2}$, then the above inequality yields
$$
B\leq N(d)B^{ d/ p_{0} } \|h\|_{p_{0},q_{0}},
\quad B^{(p_{0}- d)/p_{0}}
\leq N(d) \|h\|_{p_{0},q_{0}}
$$
and we obtain \eqref{5.6.5} again. The theorem is proved.

\begin{remark}
                                                   \label{remark 5.7.1}
In the case of $q=\infty,p=d$ an estimate of $B$
in terms of $\|h\|_{p,q}$ is given in Theorem 5.2 of \cite{Kr_86}
if $\gamma$ is the first exit time of $x_{t}$ from a ball
and in Theorem 2.17 of \cite{Kr_19} if $A_{t}=t$ and   
$c_{t}=\lambda \tr a_{t}$,
where $\lambda>0$ is a number (and $\gamma=\infty$).

\end{remark}

\begin{remark}
                                  \label{remark 5.7.2}
As in \cite{Na_15} we note that estimate
\eqref{5.6.4} also, obviously, holds 
if 
$$
|b_{t}|\leq \sum_{k=1}^{n}\kappa^{k}_{t}h_{k}(\tau_{t}
,x_{t}),
$$
where $\kappa^{k}_{t}=r_{t}^{1/q_{k}}
(\det a_{t})^{1/p_{k}}c^{\theta_{k}}_{t}$, $p_{k}\in[1,\infty],
q_{k}\in[1,\infty)$, $\theta_{k}=1-d/p_{k}-1/q_{k}\geq0$,
and $h_{k}$ are nonnegative Borel functions. In that
case the constant $C$ depends only on $d,p,q,p_{k},q_{k}$,
  $\|h_{k}\|_{p_{k},q_{k}}$, $k=1,..,n$, in a somewhat complicated way.

\end{remark}
\begin{remark}
                                  \label{remark 5.8.2}
The main case of applications of Theorem 
\ref{theorem 5.5.2} in this article is
when $p=p_{0}<\infty,q=q_{0}<\infty,\theta=\theta_{0}=0$,
$\gamma=T$, where $T$ is a fixed number, $r_{t}=1$,
$c_{t}=0$, $A_{t}= t\wedge T $,  
$$
|b_{t}|\leq (\det a_{t})^{1/p}h(t,x_{t})I_{t\leq T}.
$$
In that case $2p/(p-d)=2q$ and estimate \eqref{5.6.4} becomes
$$
E\int_{0}^{T}(\det a_{t})^{1/p}f(t,x_{t})\,dt\leq
N(d,p)\Big(T+\|hI_{(0,T)}\|_{p,q}^{ 2q}\Big)^{d/(2p)}
\|f\|_{p,q}.
$$

We finish the section with somewhat unrelated result
which we use later in Section \ref{section 5.14.1}
and
which would be a simple consequence of Theorem 4.5.1
of \cite{SV_79} if we assumed that $b$   is bounded.

\begin{lemma}
                                            \label{lemma 12.6.1}
Let $x_{t}$, $t\geq0$, be an $\bR^{d}$-valued process
on a probability space $(\Omega,\cF,P)$. Define
$\cF_{t}$ as the completion of the $\sigma$-field generated
by $x_{s}$, $s\leq t$.
Let $\sigma_{t}$ be an $\bS_{\delta}$-valued
 and $b$  be an $\bR^{d}$-valued processes which are
progressively measurable with respect to $\{\cF_{t}\}$.
Suppose that   for any $T\in(0,\infty)$
$$
\int_{0}^{T}|b_{t}|\,dt<\infty
$$
(a.s.), and for any
$C^{\infty}_{0}(\bR^{d+1})$-function $u(t,x)$
  the process
\begin{equation}
                                               \label{11.6.2}
u(t,x_{t})-\int_{0}^{t}L_{s}u(s,x_{s})\,ds
\end{equation}
is a local martingale with respect  to $\{\cF_{t}\}$,
where for $a=\sigma^{2}$
$$
L_{t}u(t,x)=\partial_{t}u(t,x)+(1/2)a^{ij}_{t}D_{ij}u(t,x)
+b^{i}_{t}D_{i}u(t,x)
$$
Then there exists a $d$-dimensional Wiener process
$(w_{t}, \cF_{t})$, $t\geq0$, 
 such that 
$$
x_{t}=x_{0}+\int_{0}^{t}\sigma_{s}\,dw_{s}+\int_{0}^{t}b_{s}\,ds.
$$
\end{lemma}

Proof. First observe that by using cut-off functions
one easily shows that \eqref{11.6.2} is a  local martingale
for any twice continuously differentiable function $u$.
Then,
we claim that the following processes are local martingales
$$
X_{t}:=x_{t}-\int_{0}^{t}b _{s} \,ds,
$$
$$
B_{t}:=x_{t}x^{*}_{t}-\int_{0}^{t}\big(a _{s} +b _{s} x^{*}_{s}
+x_{s}b^{*} _{s} \big)\,ds,
$$
$$
A_{t}:=X_{t}X_{t}^{*}-\int_{0}^{t} a _{s} \,ds.
$$

Indeed, the first two processes are obtained from 
\eqref{11.6.2} for $u=x,xx^{*}$. Concerning the last one
introduce $\gamma_{R}$ as the minimum of $\tau_{R}=\inf\{t\geq0:|x_{t}|
\geq R\}$   
and 
$$
\inf\{t\geq0:\int_{0}^{t}|b _{s} |\,ds+|B_{t}|\geq R\}.
$$
Also let
$$
\Phi_{t}=\int_{0}^{t}b _{s} I_{s<\gamma_{R}}\,ds.
$$
Observe that $X_{t\wedge\gamma_{R}}$ and $\Phi_{t}$
are bounded and simple manipulations yield
$$
A_{t\wedge\gamma_{R}}=\int_{0}^{t}X_{s\wedge\gamma_{R}}\,d
\Phi^{*}_{s}-X_{t\wedge\gamma_{R}} 
\Phi^{*}_{t}+
\int_{0}^{t}\big(d
\Phi_{s}\big)X^{*}_{s\wedge\gamma_{R}}- 
\Phi_{t}X^{*}_{t\wedge\gamma_{R}}+B_{t\wedge\gamma_{R}},
$$
which by the Lemma from Appendix 2 of \cite{Kr_77}
shows that $A_{t\wedge\gamma_{R}}$ is a martingale.

By the above claim the quadratic variation
process of the local martingale $X_{t}$
is
$$
\int_{0}^{t}a_{s} \,ds.
$$
After that our assertion   
follows directly from
  Theorem III.10.8 of \cite{Kr_95}. The lemma is proved. 

\end{remark}

\mysection{Proof of Theorem \protect\ref{theorem 9.6.4}}
                                                 \label{section 5.13.1}

Introduce  
$$
B(t)=  \|bI_{(-\infty,t)}\|_{p ,q }
 ^{ q}.
$$
\begin{lemma}
                                               \label{lemma 4.25.1}
Suppose that $p\geq q$ and
let $x_{t}$ be a solution of \eqref{11.29.2}. Then
for $0\leq s<t<s+1<\infty$ and $n=1,2,...$, we have
\begin{equation}
                                                  \label{4.25.3}
E|x_{t}-x_{s}|^{n}\leq N\big(t-s+B^{2}(t_{0}+t)-B^{2}(t_{0}+s)\big)^{nd/(2p)},
\end{equation}
where $N=N(n,d,\delta,p,\|b\| )$.
\end{lemma}
 
Proof. We may assume that $t_{0}=0$.
Then observe that for any integer $n=1,2,...$
$$
I_{n+1}:=E\Big(\int_{s}^{t}b (u,x_{u}) \,du\Big)^{n+1}
$$
$$
=(n+1)!E\int_{s\leq u_{1}\leq...\leq u_{n}}b (u_{1},x_{u_{1}})\cdot...\cdot
b (u_{n},x_{u_{n}})
$$
$$
\times E\Big(\int_{u_{n}}^{t}b (u,x_{u}) \,du\mid \cF_{u_{n}}\Big)
du_{1}...du_{n},
$$
where the conditional expectation we can estimate by using Remark
\ref{remark 5.8.2}.  

Then we get
$$
I_{n+1}\leq N(n+1)I_{n} 
\Big(t-s+\|bI_{(s,t)}\|_{p,q}^{ 2q}\Big)^{d/(2p)}
\|b \|_{p,q} ,
$$
where $N$ depends only on $d$, $p$, and $\delta$.  Here
$$
\|bI_{(s,t)}\|_{p,q}^{ 2q}=
\Big(B(t)-B(s)\Big)^{2 }\leq\ 
 B^{2 }(t)-B^{2}(s).
$$
Therefore,
$$
I_{n+1}\leq N(n+1)I_{n}\Big(t-s+B^{2}(t)-B^{2}(s)
\Big) ^{d/(2p)} \|b \|_{p,q}.
$$

The induction on $n$ yields
$$
I_{n}\leq N^{n}n!\Big(t-s+B^{2}(t)-B^{2}(s)
\Big) ^{nd/(2p)}\|b \|^{n}_{p,q}.
$$
Also, as is well known,
$$
E\Big|\int_{s}^{t}\sigma (u,x_{u}) \,dw_{u}\Big|^{ n}
\leq N(n,\delta)(t-s)^{n/2}.
$$
It follows that the left-hand side of \eqref{4.25.3}
is less than a constant $N$ times
$$
(t-s)^{n/2}+\Big(t-s+B^{2}(t)-B^{2}(s)
\Big) ^{nd/(2p)},
$$
which less than twice the  factor of $N$
in \eqref{4.25.3} because $p>d$ and $ t-s \leq 1$.
This proves the lemma.

\begin{lemma}
                                              \label{lemma 4.25.3}
Under the assumptions in Theorem \ref{theorem 9.6.4} (ii)
the set of distributions of $ x^{n }_{\cdot} $ on 
$C([0,\infty),\bR^{d })$ is tight if $p\geq q$.  

\end{lemma}

Proof. Define
$$
B_{n}(t)=\|b^{n}I_{(-\infty,t ^{n}+t)}\|^{q}_{p,q}
$$
and let $\phi^{n}(s)$ be the inverse function of $ 
\psi^{n}(t):=t ^{n}+t+B_{n}^{2}(t ^{n}+t)$.
By Lemma \ref{lemma 4.25.1} and Kolmogorov's criteria
the set of distributions of $y^{n}_{\cdot}:=
 x^{n}_{\phi^{n}(\cdot)} $ on 
$C([0,\infty),\bR^{d })$ is tight.

Observe that, as $n\to\infty$, $\psi^{n}(t)$
converges to $t_{0} +t+B^{2} (t_{0} +t)$ which is continuous and monotone.
By Polya's theorem the convergence is uniform
on any finite time interval, and hence, the
functions $\psi^{n}(t)$ are  
 equi-continuous
on any finite time interval. Now 
define

$$
\Phi(s)=\inf_{n\geq 1}\phi^{n}(s)
$$
and take $S\in(0,\infty)$. By tightness,
for any  
  $\varepsilon>0$ there is a compact set $K_{\varepsilon}$
in $C([0,S],\bR^{d })$ such that $P^{n}(\{y^{n}_{s},s\leq S\}
\in K_{\varepsilon})\geq 1-\varepsilon$ for all $n$.
Due to the uniform continuity of $\psi^{n}$
and of the elements of $K_{\varepsilon}$, the elements of
$$
\hat K_{\varepsilon}
:=\{\{f(\psi^{n}(t)),t\leq \Phi(S)\}: \{f(s),s\leq S\}
\in K_{\varepsilon},n=1,2,...\}
$$
are uniformly continuous and, of course, uniformly bounded,
so that $\hat K_{\varepsilon}$ is a compact set 
in $C([0,\Phi(S)],\bR^{d})$
and
$$
P(\{y^{n}_{\psi^{n}(t)},t\leq \Phi(S)\}\in\hat K_{\varepsilon})
\geq 1-\varepsilon.
$$
It only remains to observe that
 $y^{n}_{\psi^{n}(t)}=x^{n}_{t}$, $S$ is arbitrary,
and $\Phi(S)\to\infty$ as $S\to\infty$. The lemma is proved.

This takes care of part of assertion (ii) of
Theorem \ref{theorem 9.6.4}. To deal with the rest 
 we rely on the following
  results  due to A. V. Skorokhod
(see Ch.~1, \S6 and Ch.~2, \S3 in \cite{Sk_61}).  
\begin{lemma}
                                               \label{lemma 4.23.1}
  Suppose that $d_1$-dimensional random
 processes $\xi^{n}_{t}$  $(t\geq 0, n =  
1,2, . . .)$ are defined
 on some probability spaces. Assume that for each $T> 0$ and
$\varepsilon >  0 $
\begin{equation}
                                                  \label{5.10.5}
\lim_{c\to\infty} \sup_{n} \sup_{t\leq T} P^{n} (|\xi^{n}_{t}|>c) = 0,
\end{equation}
\begin{equation}
                                                  \label{5.10.6}
\lim_{h\downarrow 0} \sup_{n} \sup_{\substack{t_{1},t_{2}\leq T\\
|t_{1}-t_{2}|\leq h}} P^{n}(|\xi^{n}_{t_{1}}-\xi^{n}_{t_{2}}|>
\varepsilon)=0.
\end{equation}
Then, one can choose a sequence of numbers $n'\to\infty$,
 a probability space, and
random processes $\tilde \xi_{t},\tilde \xi_{t}^{n'}$
 defined on this probability space such that all finite-dimensional
distributions of $\tilde \xi_{t}^{n'}$ coincide with 
the corresponding finite-dimensional
distributions of $\xi_{t}^{n'}$ and 
$$
P (|\tilde \xi_{t}-\tilde \xi_{t}^{n'}|) \to 0   
$$
as $n'\to\infty$ for any $\varepsilon>0$ and $t\geq0$.
\end{lemma}

\begin{lemma}
                                               \label{lemma 4.23.2}
Suppose the assumptions of Lemma  \ref{lemma 4.23.1} 
 are satisfied
and $\xi^{n}_{t}$ are defined on the same probability space.
 Also, suppose
that   $d_{1}$-dimensional Wiener processes $(w^{n}_{t},\cF^{n}_{t})$
 are defined on this 
probability space. Assume that the functions 
$\xi^{n}_{t}(\omega)$ are bounded on $[0,\infty)\times\Omega$
uniformly in n and that the stochastic integrals 
$$
I^{n}_{t}:=\int_{0}^{t}\xi^{n}_{s}\,dw^{n}_{s}
$$
 are defined for $t\geq0$.
Finally, let 
\begin{equation}
                                                     \label{5.14.1}
 \xi^{n}_{t}\to \xi^{0}_{t},\quad w^{n}_{t}\to w^{0}_{t}
\end{equation}
in probability as $n\to\infty$
for each $t\geq0$. Then $I^{n}_{t}\to I^{0}_{t}$
in probability as $n\to\infty$
for each $t\geq0$.
\end{lemma}

\begin{remark}
                                                 \label{remark 5.11.1}
As it follows from the proof of Lemma \ref{lemma 4.23.2}
given in \cite{Sk_61}
  we need conditions
\eqref{5.10.5}, \eqref{5.10.6}, and \eqref{5.14.1} to hold
only for $t,t_{1},t_{2}$ restricted to a set
of full measure in order for the assertion of the lemma to be true.
\end{remark}

\begin{lemma}
                                             \label{lemma 5.12.1}
Let $\bR^{2d}$-valued processes $(x^{i}_{t},w^{i}_{t})$, $t\geq0$,
$i=1,2$
defined on perhaps different probability spaces have the same
finite-dimensional distributions. Define 
$\cF^{i}_{t}$ as the  completion of
  $\sigma(x^{i}_{s},w^{i}_{s}:s\leq t)$ and assume that $w^{1}_{t}$
is a Wiener process with respect to   $\cF^{1}_{t}$. Also suppose that
(a.s.) for all $t\geq 0$
\begin{equation}
                                              \label{5.13.1}
x^{1}_{t}=\int_{0}^{t}\sigma(s,x^{1}_{s})\,dw^{1}_{s}+
\int_{0}^{t}b(s,x^{1}_{s})\,ds.
\end{equation}
Then $x^{2}_{t},w^{2}_{t}$ have   modifications 
(called again $x^{2}_{t},w^{2}_{t}$) such that $w^{2}_{t}$
is a Wiener process with respect to   $\cF^{2}_{t}$
and (a.s.) for all $t\geq 0$
\begin{equation}
                                              \label{5.13.2}
x^{2}_{t}=\int_{0}^{t}\sigma(s,x^{2}_{s})\,dw^{2}_{s}+
\int_{0}^{t}b(s,x^{2}_{s})\,ds.
\end{equation}
\end{lemma}

Proof. Fix $T\in(0,\infty)$ and $\varepsilon\in(0,1)$.
Since the trajectories of $(x^{1}_{t},w^{1}_{t})$
are continuous, there exists a compact set $K
\subset C([0,T],\bR^{2d})$ such that
$$
P((x^{1}_{\cdot\wedge T},w^{1}_{\cdot\wedge T})\in K)\geq 1-\varepsilon.
$$
Hence, there is a constant $N$ and a continuous
function $w(t)$, $t\in[0,T]$, such that $w(0)=0$ and
with probability larger than $1-\varepsilon$ for any $s,t\in[0,T]$
\begin{equation}
                                              \label{5.12.2}
|(x^{1}_{s},w^{1}_{s})|\leq N,\quad
|(x^{1}_{s},w^{1}_{s})-(x^{1}_{t},w^{1}_{t})|\leq w(|t-s|).
\end{equation}

It follows that \eqref{5.12.2} holds for rational $s,t$
if we replace $(x^{1},w^{1})$ with $(x^{2},w^{2})$.
Then by continuity $(x^{2}_{t},w^{2}_{t})$ is extended
to all $t\in[0,T]$. The extensions coincide
with the original ones (a.s.) for any $t$ because
of the stochastic continuity of the original $(x^{2}_{t},w^{2}_{t})$.
This is done on events   whose  probabilities tend
to one. Because of the arbitrariness of $T$ we may assume that
$(x^{2}_{t},w^{2}_{t})$ is continuous in $t$ with probability one.

By Remark  \ref{remark 5.8.2}  and by the coincidence
of finite dimensional distributions (and by the measurability of
$x^{2}_{t}$ due to its continuity)
for any $T\in[0,\infty)$, Borel $f(t,x)\geq0$,  
\begin{equation}
                                                  \label{5.13.3}
E\int_{0}^{T}
f(t,x^{2}_{t})\,dt\leq N\|fI_{(0,T)}\|_{p,q},
\end{equation}
where $N$ is independent of $f$.

Furthermore, if $\alpha(t,x)$ is a continuous
$d\times d$ symmetric matrix-valued, $\beta(t,x)$ is 
a continuous $\bR^{d}$-valued, then the distributions of
$$
(x^{i}_{t},\int_{0}^{t}\alpha(s,x^{i}_{s})\,dw^{i}_{s},\int_{0}^{t}
\beta(s,x^{i}_{s})\,ds),\quad i=1,2,
$$
coincide, because the integrals can be approximated by
integral sums. This coincidence also holds for $\alpha=\sigma$
and $\beta=b$ due to \eqref{5.13.3} and the possibility 
of approximation. Hence for each $t$ with probability
one \eqref{5.13.2} holds due to \eqref{5.13.1}.
But then with probability one it holds for all $t$,
because both sides of \eqref{5.13.2} are continuous.
The lemma is proved.

{\bf Proof of Theorem  \ref{theorem 9.6.4}}.
Due to the possibility to use
 mollifiers we see that assertion (ii)
implies (i).
In the proof of (ii), thanks to   
 Lemma \ref{lemma 4.25.3}, we need only prove
the assertion concerning
the convergence of finite dimensional distributions. 

Having in mind Lemma \ref{lemma 4.23.1} define for $M>0$
$$
\xi^{n}_{t}=\int_{0}^{t}b^{n}(t^{n}+s,x^{n}_{s})\,ds,\quad
\xi^{nM}_{t}=\int_{0}^{t}b^{n}(t^{n}+s,x^{n}_{s})
I_{|b^{n}( t^{n}+s, x^{n}_{s})|\leq M}\,ds.
$$

Since the derivative of $\xi^{nM}_{t}$ is bounded,
both conditions \eqref{5.10.5} and \eqref{5.10.6}
are satisfied for $\xi^{nM}_{t}$. Furthermore,
$$
P^{n}\Big(\int_{0}^{T}|b^{n}(t^{n}+s,x^{n}_{s})|
I_{|b^{n}(t^{n}+s,x^{n}_{s})|\geq M}\,ds>\varepsilon\Big)\leq \varepsilon^{-1}
N\|b^{n}I_{|b^{n}|\geq M}\|_{p,q},
$$
where $N$ is independent of $n$ and $\varepsilon$. Since $b^{n}\to b$
in the $\|\cdot\|_{p,q}$-norm, the latter quantity can be made
as small as we like on the account of choosing $M$
 large enough. Therefore, Lemma \ref{lemma 4.23.1} is applicable
to $\xi^{n}_{t}$. It is, obviously, also applicable to
$$
\eta^{n}_{t}=x^{n}+\int_{0}^{t}\sigma^{n}(t^{n}+s,x^{n}_{s}) \,dw^{n}_{s}.
$$
Hence, there is a subsequence, which by common abuse of notation
we identify with the original one, a probability space and
random $\bR^{2d}$-valued
processes $(\tilde x^{n}_{t},\tilde w^{n}_{t})$,
$(\tilde x^{0}_{t},\tilde w^{0}_{t})$  
defined on this probability space such that all finite-dimensional
distributions of $(\tilde x^{n}_{t},\tilde w^{n}_{t})$ coincide with 
the corresponding finite-dimensional
distributions of $(x_{t}^{n },w^{n}_{t})$ and 
\begin{equation}
                                                       \label{5.13.6}
P (|(\tilde x^{n}_{t},\tilde w^{n}_{t})
-(\tilde x^{0}_{t},\tilde w^{0}_{t})|\geq \varepsilon) \to 0   
\end{equation}
as $n \to\infty$ for any $\varepsilon>0$ and $t\geq0$. 
Furthermore, for any $T\in(0,\infty)$ there exists
a continuous function $w(t)$, $t\in[0,T]$, such that
$w(0)=0$ and for all $n\geq 0$, $s,t\leq T$,
\begin{equation}
                                                       \label{5.13.4}
E|\phi(\tilde x^{n}_{t})-\phi(\tilde x^{n}_{s})|\leq w(|t-s|),
\end{equation}
where $\phi(x)=x/(1+|x|)$.

For $n\geq0$ introduce $\tilde \cF^{n}_{t}$ as the completion
of $\sigma(\tilde x^{n}_{s}, \tilde w^{n}_{s},s\leq t)$. It is easy to see,
using Kolmogorov's continuity criterion, that
$\tilde w^{0}_{t}$ admits a continuous modification
$\hat w^{0}_{t}$ such that
$\{\hat w^{0}_{t},\tilde \cF^{0}_{t}\}$  is a Wiener
process.

By Lemma \ref{lemma 5.12.1}, for each $n\geq1$, the process
$(\tilde x^{n}_{t},\tilde w_{t}^{n})$
admits a continuous
modification denoted by   $(\hat x^{n}_{t},\hat w_{t}^{n})$
such  that 
 $(\hat w_{t}^{n},
\tilde \cF_{t}^{n})$ is a Wiener process and  (a.s) for all $t\geq0$
\begin{equation}
                                                     \label{5.12.1}
\hat x^{n}_{t}=x_{n}+\int_{0}^{t}\sigma^{n}(t_{n}+s,
\hat x^{n}_{s})\,d\hat w^{n}_{s}+\int_{0}^{t}b^{n}(t_{n}+s,
\hat x^{n}_{s})\,ds.
\end{equation}

In light of \eqref{5.13.6} and \eqref{5.13.4} we have
\begin{equation}
                                                       \label{5.17.1}
P (|(\hat x^{n}_{t},\hat w^{n}_{t})
-(\tilde x^{0}_{t},\tilde w^{0}_{t})|\geq \varepsilon) \to 0   
\end{equation}
as $n \to\infty$ for any $\varepsilon>0$ and $t\geq0$ and 
for all $n\geq 1$, $s,t\leq T$,
\begin{equation}
                                                       \label{5.17.2}
E|\phi(\hat x^{n}_{t})-\phi(\hat x^{n}_{s})|\leq w(|t-s|).
\end{equation}

Now the fact that $\tilde x^{0}_{t}$ may be not measurable
in $t$ causes some problems. However, observe that, owing to \eqref{5.17.1},
 $\phi(\hat x^{n}_{t})$ form a Cauchy sequence
in $L_{1}(\Omega\times[0,T])$ and, hence, converges in that space
to $\phi(\hat x^{0}_{t})$, where $\hat x^{0}_{t}$ is measurable
with respect to $(\omega,t)$. By Fubini's theorem
there is a set $\cS\subset [0,\infty)$ of full measure
such that, for any $t\in\cS$, $\hat x^{0}_{t}=\tilde x^{0}_{t}$
(a.s.).  Without losing the above properties we set
$\hat x^{0}_{t}=0$ for $t\not\in \cS$ and 
then, for any $s,t\geq 0$, $\hat w^{0}_{t+s}-
\hat w^{0}_{t}$ is indepenent of 
$(\hat x^{0}_{r},\hat w^{0}_{r}),r\leq t$. 

Now we note that  \eqref{5.17.2} remains 
valid for $n=0$ and \eqref{5.17.1}  remains 
valid if we replace $(\tilde x^{0}_{t},\tilde w^{0}_{t})$
by $(\hat x^{0}_{t},\hat w^{0}_{t})$ and 
restrict the ranges of $t,s$ to $t,s\in\cS$.
This is done to accommodate Remark \ref{remark 5.11.1}.
 Then by Lemma
\ref{lemma 4.23.2} for any  $t\geq0$ and  continuous
$d\times d$ symmetric matrix-valued
 $\alpha(t,x)$ we have 
\begin{equation}
                                                     \label{5.13.7}
 \int_{0}^{t}\alpha(s,\hat x_{s}^{n})\,d\hat w^{n}_{s}
\to \int_{0}^{t}\alpha(s,\hat x_{s}^{0})\,d\hat w^{0}_{s}
\end{equation}
as $n\to\infty$ in probability.
We want to use this  to pass to the limit in the stochastic term
in \eqref{5.12.1}.
But first observe that
by Remark  \ref{remark 5.8.2}  
for any $T\in[0,\infty)$, Borel $f(t,x)\geq0$, and $n\geq1$
\begin{equation}
                                                  \label{4.26.1}
E\int_{0}^{T}
f(t,\hat x^{n}_{t})\,dt\leq N\|fI_{(0,T)}\|_{p,q},
\end{equation}
where $N$ is independent of $f$ and $n$. The convergence
in probability implies that \eqref{4.26.1} holds for
$n=0$ as well with the same constant $N$, first for 
nonnegative $f\in C^{\infty}_{0}
(\bR^{d+1})$ and then, due to general measure-theoretic arguments,
for any Borel nonnegative $f$.

We claim that on the account of \eqref{4.26.1},
if Borel functions $g^{n}$ converge to $g$ in 
the $\|\cdot\|_{p,q}$-norm,
then
\begin{equation}
                                                  \label{4.26.2}
E\int_{0}^{T}|
g^{n}(t,\hat x^{n}_{t})-
g (t,\hat x^{0}_{t})|\,dt\to 0.
\end{equation}
To prove \eqref{4.26.2}    take
$\varepsilon>0$ and $g_{\varepsilon}\in C^{\infty}_{0}
(\bR^{d+1})$ such that
$$
\|g-g_{\varepsilon}\|_{p,q}\leq \varepsilon.
$$
For $g_{\varepsilon}$ in place of $g$, \eqref{4.26.2}
follows from the convergence in probability
of $\hat x^{n}_{t}$ to $\hat x^{0}_{t}$ for $t\in\cS$. 
After that it only remains to observe that the limit
of the error of the substitution in
  \eqref{4.26.2} is less than $2N\varepsilon$ 
owing to \eqref{4.26.1}. It follows, in particular, that in probability
\begin{equation}
                                                  \label{5.13.5}
\sup_{t\leq T}\Big|\int_{0}^{t}b^{n}(t_{n}+s,
\hat x^{n}_{s})\,ds
-\int_{0}^{t}b (t_{0}+s,
\hat x^{0}_{s})\,ds\Big|\to 0.
\end{equation}

Coming back to the stochastic part note that for any $t\geq0$ and $ 
c\in(0,\infty)$ 
$$
\nlimsup_{n\to\infty}E\Big|\int_{0}^{t}\sigma^{n}(t_{n}+s,
\hat x^{n}_{s})\,d\hat w^{n}_{s}-
\int_{0}^{t}\alpha(s,\hat x_{s}^{n})\,d\hat w^{n}_{s}\Big|^{2}
$$
$$
=\nlimsup_{n\to\infty}E\int_{0}^{t}\|\sigma^{n}(t_{n}+s,
\hat x^{n}_{s})-\alpha(s,\hat x_{s}^{n})\|^{2}\,ds
$$
$$
\leq N \sup_{n}\int_{0}^{t}P(|\hat x^{n}_{s}|>c)\,ds+
N\lim_{n\to\infty}\big\|\big(\sigma^{n}(t_{n}+\cdot,
\cdot)-\alpha(\cdot,\cdot)\big)I_{[0,t]\times B_{c}}\big\|_{p,q} 
$$
$$
=N \sup_{n}\int_{0}^{t}P(|\hat x^{n}_{s}|>c)\,ds+
N \big\|\big(\sigma (t_{0}+\cdot,
\cdot)-\alpha(\cdot,\cdot)\big)I_{[0,t]\times B_{c}}\big\|_{p,q},
$$
where the constants $N$ are independent of $t$ and $c$.
The last quantity also dominates
$$
E\Big|\int_{0}^{t}\sigma (t_{0}+s,
\hat x^{0}_{s})\,d\hat w^{0}_{s}-
\int_{0}^{t}\alpha(s,\hat x_{s}^{0})\,d\hat w^{0}_{s}\Big|^{2}.
$$
This and \eqref{5.13.7} show how, for any given $\varepsilon,\delta>0$,
 to choose $c$ and a continuous $\alpha$  in order to have that
$$
\nlimsup_{n\to\infty}P\Big(
\Big|\int_{0}^{t}\sigma^{n}(t_{n}+s,
\hat x^{n}_{s})\,d\hat w^{n}_{s}-
\int_{0}^{t}\sigma (t_{n}+s,
\hat x^{0}_{s})\,d\hat w^{0}_{s}\Big|>\varepsilon  
\Big)\leq \delta.
$$

Upon combining this with \eqref{5.13.5} and coming back to
\eqref{5.12.1} we conclude that for any $t$ (a.s.)
$$
\tilde x^{0}_{t}=x_{0}+\int_{0}^{t}\sigma (t_{0}+s,
\hat x^{0}_{s})\,d\hat w^{0}_{s}+\int_{0}^{t}b (t_{0}+s,
\hat x^{0}_{s})\,ds=:y_{t}.
$$
In particular, this means that $\tilde x^{0}_{t}$
admits a continuous modification $y_{t}$. In turn, it
allows us to replace in the above equation
$\hat x^{0}_{s}$ with $y_{t}$, because for any $s\in \cS$,
$\hat x^{0}_{s}=\tilde x^{0}_{s}=y_{s}$ (a.s.) and therefore
 $\hat x^{0}_{s} =y_{s}$ for almost all $(\omega,s)$.
This, of course, brings the proof of the theorem to an end.

\mysection{Markov processes corresponding to $\sigma,b$}
                                        \label{section 5.14.1}

We are going to use the results in \cite{Kr_73}
applied in the case when the semicompactum $E$
is $\bR^{d+1}$, that is when the $t$-variable
is considered just as one of coordinates of points $(t,x)
\in\bR^{d+1}$.

Let $\Omega$ be the set of $\bR^{d+1}$-valued
 continuous function $(t_{0}+t,x_{t})$, $t_{0}\in \bR$,
defined for $t\in[0,\infty)$
For $\omega=\{(t_{0}+t,x_{t}),t\geq0 \}$, define
$\sft_{t}(\omega)=t_{0}+t$, $x_{t}(\omega)=x_{t}$,
and set $\cN_{t}=\sigma((\sft_{s},x_{s}),s\leq t)$,
$\cN=\cN_{\infty}$. Denote by $\sfT$ the set of stopping times
relative to $\{\cN_t\}$.
In the following theorem we use the terminology from
\cite{Dy_63}.

\begin{theorem}
                                            \label{theorem 4.27.1}
 On $\bR^{d+1}$ there exists a strong Markov process
$$
X=\{(\sft_{t},x_{t}),\infty,\cN_{t}, P_{t,x})
$$
such that the process 
$$
X_{1}=\{(\sft_{t},x_{t}),\infty,
\cN_{t+}, P_{t,x})
$$
 is Markov and for any $(t,x)\in\bR^{d+1}$
there exists a $d$-dimensional Wiener process $w_{t}$, $t\geq0$,
which is a Wiener process relative to $\bar \cN_{t}$,
where $\bar \cN_{t}$ is the completion of $\cN_{t}$
with respect to $P_{t,x}$, and such that with 
$P_{t,x}$-probability one, for
all $s\geq 0$, $\sft_{s}=t+s$ and
\begin{equation}
                                            \label{4.27.1}
x_{s}=x+\int_{0}^{s}\sigma(t+u,x_{u})\,dw_{u}
+\int_{0}^{s}b(t+u,x_{u})\,du.
\end{equation}
\end{theorem}

Proof. Define $a=\sigma^{2}$,
$$
Lu(t,x)=\partial_{t}u(t,x)+(1/2)a^{ij}D_{ij}u(t,x)+b^{i}D_{i}u(t,x)
$$
and introduce $\Pi_{t,x}$ as the set of probability
measures on $(\Omega,\cN)$ such that $P((\sft_{0},x_{0})=(t,x))=1$,
\begin{equation}
                                                     \label{4.27.5}
E\int_{0}^{T}|b(\sft_{t},x_{t})|\,dt<\infty,\quad\forall T<\infty,
\end{equation}
and the process
$$
\eta_{t}(u)=u(\sft_{t},x_{t})-\int_{0}^{t}Lu(\sft_{s},x_{s})\,ds
$$
is a martingal relative to $\{\cN_{t}\}$ for all $u\in C^{\infty}
_{0}(\bR^{d+1})$.

According to Lemma \ref{lemma 12.6.1}, if $P_{t,x}\in\Pi_{t,x}$, then
the assertion of the theorem regarding \eqref{4.27.1}
holds and \eqref{4.27.5} is true.
Therefore, by Theorem 2 of \cite{Kr_73}   to prove
the present theorem, it suffices to show that $\Pi_{t,x}\ne\emptyset$
and
$\{\Pi_{t,x}\}$ is a Markov system relative to $(\sf T,\cN_t)$
and $([0,\infty),\cN_{t+})$.

That $\Pi_{t,x}\ne\emptyset$ follows from Theorem \ref{theorem 9.6.4} (i).
Let us prove that $\{\Pi_{t,x}\}$ is a $B$-system. To achieve this,
as it follows from \cite{Kr_73}, it suffices to show that
if $(t_{n},x_{n})\to (t,x)$ and $P^{n}\in\Pi_{t_{n},x_{n}}$,
then there exists a subsequence $n(k)\to\infty$ and $P^{0}\in\Pi_{t,x}$
such that for any $f\in C^{\infty}_{0}(\bR^{d+2})$
$$
E^{ n(k) }\exp\Big(\int_{0}^{\infty}e^{-t}f(t,\sft_{t},x_{t})\,dt\Big)
\to E^{0}\exp\Big(\int_{0}^{\infty}e^{-t}f(t,\sft_{t},x_{t})\,dt\Big),
$$
where $E^{ n(k) },E^{0}$ are the expectation signs with respect
to $P^{ n(k) },P^{0}$, respectively.
The reader will easily derive this property   
from Theorem \ref{theorem 9.6.4} (ii) by 
using Taylor's series and
observing that
$$
E\Big(\int_{0}^{\infty}e^{-t}f(t,\sft_{t},x_{t})\,dt\Big)^{n}
$$
$$
=E \int_{0}^{\infty}...\int_{0}^{\infty}
e^{-t_{1}}f(t_{1},\sft_{t_{1}},x_{t_{1}})\cdot...
\cdot e^{-t_{n}}f(t_{n},\sft_{t_{n}},x_{t_{n}})\,
dt_{1}\cdot...\cdot dt_{n}.
$$

What remains is to prove that for $(\sfT,\cN_{t})$
and $([0,\infty),\cN_{t+})$ the conditions 2) and 3)
are satisfied
of the definition of Markov system in \cite{Kr_73}.
This is done by almost literally repeating the corresponding part
of the proof
of Theorem 3 of \cite{Kr_73}. One need only  replace there
$x_{t}$ with $(\sft_{t},x_{t})$. The theorem is proved.

{\bf Acknowledgment}. The author is sincerely grateful
to A.I. Nazarov, who pointed out an error in the first version
of the article,  to Hongjie Dong and Doyoon Kim for spotting several
misprints bordering with errors, and to Xicheng Zhang whose comment
allowed the author to avoid an incorrect statement.

\end{document}